\documentclass[12pt, reqno]{amsart}
\usepackage{amsmath,amssymb}

\usepackage{amsmath, amsthm, amscd, amsfonts, amssymb, graphicx}
\usepackage[bookmarksnumbered,plainpages]{hyperref}

\newtheorem{theorem}{Theorem}[section]

\newtheorem{corollary}[theorem]{Corollary}
\theoremstyle{definition}

\theoremstyle{remark}

\numberwithin{equation}{section}

\begin{document}
\title[On approximate cubic homomorphisms]
{On approximate cubic homomorphisms}
\author{M. Eshaghi Gordji and M. Bavand Savadkouhi}
\address{Department of Mathematics,
 Semnan University, P. O. Box 35195-363, Semnan, Iran} \email{
madjid.eshaghi@gmail.com \& bavand.m@gmail.com}
\keywords{Cubic functional equation; Homomorphism;
 Hyer-Ulam-Rassias stability}
\subjclass[2000]{Primary 39B52, Secondary 39B82, 46H25}
\smallskip
\begin{abstract}
In this  paper, we investigate the generalized Hyers--Ulam--Rassias
stability of the system of functional equations
$$\hspace {1.5 cm}\begin{cases}
f(xy)=f(x)f(y), & \\
&\qquad\qquad.\\
f(2x+y)+f(2x-y)=2f(x+y)+2f(x-y)+12f(x), &
\end{cases}$$  on  Banach
algebras. Indeed we establish the superstability of above system by
suitable control functions.
\end{abstract}
\maketitle
\section{Introduction}
A definition of stability in the case of homomorphisms between
metric groups was suggested by a problem by S. M. Ulam [21] in 1940.
Let $(G_1,.)$ be a group and let $(G_2,\ast)$ be a metric group with
the metric $d(.,.)$. Given $\epsilon>0$, does there exist a
$\delta>0$ such that if a mapping $h:G_1\longrightarrow G_2$
satisfies the inequality $d(h(x.y),h(x)\ast h(y))<\delta$ for all
$x,y\in G_1$, then there exists a homomorphism $H:G_1\longrightarrow
G_2$ with $d(h(x),H(x))<\epsilon$ for all $x\in G_1$? In this case,
the equation of homomorphism $h(x.y)=h(x)\ast h(y)$ is called
stable. In the other hand we are looking for situations when the
homomorphisms are stable, i.e., if a mapping is an approximate
homomorphism, then there exists an exact homomorphism near it. The
concept of stability for a functional equation arises when we
replace the functional equation by an inequality which acts as a
perturbation of the equation. In 1941, D. H. Hyers [8] gave a
positive answer to the question of Ulam for Banach spaces. Let
$f:E_1\longrightarrow E_2$ be a mapping between Banach spaces such
that
$$\|f(x+y)-f(x)-f(y)\|\leq \delta$$ for
all $x,y\in E_1$ and for some $\delta\geq0$. Then there exists a
unique additive mapping $T:E_1\longrightarrow E_2$ satisfying
$$\|f(x)-T(x)\|\leq \delta$$ for all $x\in E_1$. Moreover,
if $f(tx)$ is continuous in $t$ for each fixed $x\in E_1$, then the
mapping $T$ is linear.  Th. M. Rassias [20] succeeded in extending
the result of Hyers' Theorem by weakening the condition for the
Cauchy difference controlled by $(\|x\|^p+\|y\|^p)$, $p\in [0,1)$ to
be unbounded.   This condition has been assumed further till now,
through the complete Hyers direct method, in order to prove
linearity for generalized Hyers-Ulam stability problem forms. A
number of mathematicians were attracted to the pertinent stability
results of  Th. M. Rassias [20],  and stimulated to investigate the
stability problems of functional equations. The stability phenomenon
that was introduced and proved by  Th. M. Rassias is called
Hyers--Ulam--Rassias stability. And then the stability problems of
several functional equations have been extensively investigated by a
number of authors and there are
many interesting results concerning this problem (see [5-7], [9], [12-13] and [16-18]).\\
D.G.
Bourgin [4] is the first mathematician dealing with stability of
(ring) homomorphism $f(xy)=f(x)f(y)$. The topic of approximate
homomorphisms was studied by a number of mathematicians, see
[2,3,10,14,15,19] and references therein. \\
Jun and Kim [11] introduced the following functional equation
$$f(2x+y)+f(2x-y)=2f(x+y)+2f(x-y)+12f(x)\eqno(1.1)$$
and they established the general solution and generalized
Hyers--Ulam--Rassias stability problem for this functional equation.
It is easy to see that the function $f(x)=c x^3$ is a solution of
the functional equation $(1.1).$ Thus, it is natural that $(1.1)$ is
called a cubic functional equation and every solution of the cubic
functional equation is said to be a cubic function.

Let $R$ be a ring. Then a mapping $f:R\longrightarrow R$ is called a
cubic homomorphism if $f$ is a cubic function satisfying
$$f(ab)=f(a)f(b),\eqno (1.2)$$ for all $a,b\in R.$ For instance, let
$R$ be commutative, then the mapping $f:R\longrightarrow R$ defined
by $f(a)=a^3 (a\in R),$ is a cubic homomorphism. It is easy to see
that a cubic homomorphism is a ring homomorphism  if and only if it
is zero function. In this paper we study the stability of cubic
homomorphisms  on Banach algebras.\\

\section{Main results }
In the following we suppose that $A$ is a  normed algebra, $B$ is a
Banach algebra and $f$ is a mapping from $A$ into $B$, and
$\varphi,\varphi_1,\varphi_2$ are maps from $A\times A$ into
$\mathbb{R}^+$. Also, we put $0^p=0$ for $p\leq 0.$
\begin{theorem}
Let
$$\|f(xy)-f(x)f(y)\|\leq \varphi_1(x,y), \eqno(2.1)$$
and
$$\|f(2x+y)+f(2x-y)-2f(x+y)-2f(x-y)-12f(x)\|\leq
\varphi_2(x,y), \eqno(2.2)$$ for all $x,y\in A.$ Assume that the
series
$$\Psi(x,y)=\sum_{i=0}^\infty
\frac{\varphi_2(2^ix,2^iy)}{2^{3i}}$$ converges and that
$$\lim_{n\longrightarrow
\infty}\frac{\varphi_1(2^nx,2^ny)}{2^{6n}}=0,$$ for all $x,y\in
A$. Then there exists a unique cubic homomorphism $~T:A
\longrightarrow A$ such that
$$\|T(x)-f(x)\|\leq \frac{1}{16}\Psi(x,0), \eqno(2.3)$$ for all
$x\in A$.
\end{theorem}
\begin{proof}
Setting $y=0$ in $(2.2)$ yields
$$\|2 f(2x)-2^4 f(x)\|\leq \varphi_2(x,0)~, \eqno (2.4)$$
and then dividing by $2^4$ in $(2.4)$ to obtain
$$ \|\frac{f(2x)}{2^3}-f(x)\| \leq
\frac{\varphi_2(x,0)}{2.2^3}~, \eqno(2.5)$$ for all $x\in A$. Now
by induction we have
$$\|\frac{f(2^nx)}{2^{3n}}-f(x)\|\leq \frac{1}{2.2^3} \sum_{i=0}^{n-1}
\frac{\varphi_2(2^ix,0)}{2^{3i}}~. \eqno(2.6)$$ In order to show
that the functions $T_n(x)=\frac{f(2^nx)}{2^{3n}}$ is a convergent
sequence, we use the Cauchy convergence criterion. Indeed, replace
$x$ by $2^mx$ and divide  by $2^{3m}$ in $(2.6)$, where $m$ is an
arbitrary positive integer. We find that
$$\|\frac{f(2^{n+m}x)}{2^{3(n+m)}}-\frac{f(2^mx)}{2^{3m}}\| \leq
\frac{1}{2.2^3}\sum_{i=0}^{n-1}
\frac{\varphi_2(2^{i+m}x,0)}{2^{3(i+m)}}=\frac{1}{2.2^3}\sum_{i=m}^{n+m-1}
\frac{\varphi_2(2^ix,0)}{2^{3i}}$$ for all positive integers m,n.
Hence by the Cauchy criterion the limit $T(x)=\lim_{n \to \infty}
T_n(x)$ exists for each $x\in A$. By taking the limit as
$~n\longrightarrow\infty$ in $(2.6)$, we see that $\|T(x)-f(x)\|\leq
\frac{1}{2.2^3}\sum_{i=0}^\infty
\frac{\varphi_2(2^ix,0)}{2^{3i}}=\frac{1}{16}\Psi(x,0)$ and $(2.3)$
holds for all $x\in A$. If we replace $x$ by $2^nx$ and $y$ by
$2^ny$ respectively  in $(2.2)$ and divide by $2^{3n}$, we see that
\begin{align*}
&\|\frac{f(2.(2^nx)+2^ny)}{2^{3n}}+\frac{f(2.(2^nx)-2^ny)}{2^{3n}}-2\frac{f(2^nx+2^ny)}{2^{3n}}-2\frac{f(2^nx-2^ny)}{2^{3n}}\\
&-12\frac{f(2^nx)}{2^{3n}}\|\leq
\frac{\varphi_2(2^nx,2^ny)}{2^{3n}}~.
\end{align*}
 Taking the limit as
$n\longrightarrow \infty$, we find that $T$ satisfies $(1.1)$ (see
Theorem 3.1 of [11]). On the other hand we have
\begin{align*}
&\|T(xy)-T(x).T(y)\|=\|\lim_{n \to
\infty}\frac{f(2^nxy)}{2^{3n}}-\lim_{n \to
\infty}\frac{f(2^nx)}{2^{3n}}.\lim_{n \to
\infty}\frac{f(2^ny)}{2^{3n}}\|\\
&=\lim_{n \to
\infty}\|\frac{f(2^nx2^ny)}{2^{6n}}-\frac{f(2^ny)f(2^ny)}{2^{6n}}\|\\
& \leq \lim_{n \to \infty}\frac{\varphi_1(2^nx,2^ny)}{2^{6n}}=0~.
\end{align*}
for all $x,y \in A.$ We find that $T$ satisfies $(1.2)$. To prove
the uniqueness property of $T$, let $\acute{T}:A\rightarrow A$ be a
functions satisfies
$\acute{T}(2x+y)+\acute{T}(2x-y)=2\acute{T}(x+y)+2\acute{T}(x-y)+12\acute{T}(x)$
 and $\|\acute{T}(x)-f(x)\|\leq \frac{1}{16}\Psi(x,0) ~.$ Since $T, \acute{T}$ are cubic, then we have
$$T(2^nx)=2^{3n}T(x), \acute{T}(2^nx)=2^{3n}\acute{T}(x)$$
for all $x \in A$, hence,
\begin{align*}
\|T(x)-\acute{T}(x)\|&=\frac{1}{2^{3n}}\|T(2^n x)-\acute{T}(2^n x)\| \\
&\leq \frac{1}{2^{3n}}(\|T(2^n x)-f(2^n x)\|+\|\acute{T}(2^n x)-f(2^n x)\|) \\
&\leq \frac{1}{2^{3n}} (\frac{1}{2.2^3}\Psi(2^n x,0)+\frac{1}{2.2^3}\Psi(2^n x,0))\\
&=\frac{1}{2^{3(n+1)}} \Psi(2^n x,0)=\frac{1}{2^{3(n+1)}} \sum_{i=0}^{\infty} \frac{1}{2^{3i}} \varphi_2(2^{i+n}x,0) \\
&=\frac{1}{2^3} \sum_{i=0}^{\infty}\frac{1}{2^{3(i+n)}}
\varphi_2(2^{i+n}x,0)=\frac{1}{2^3}
\sum_{i=n}^{\infty}\frac{1}{2^{3i}}\varphi_2(2^i x,0)~.
\end{align*}
By taking  $~n\rightarrow \infty$ we get, $T(x)=\acute{T}(x).$
\end{proof}

\begin{corollary}
Let $\theta_1$ and $\theta_2$ be nonnegative real numbers, and let
$p\in (-\infty,3)$. Suppose that
$$\|f(xy)-f(x)f(y)\|\leq \theta_1 ~,$$
$$\|f(2x+y)+f(2x-y)-2f(x+y)-2f(x-y)-12f(x)\|\leq \theta_2(
\|x\|^p+\|y\|^p)~,$$ for all $x,y\in A$. Then there exists a
unique  cubic homomorphism $~T:A \longrightarrow A$ such that
$$\|T(x)-f(x)\|\leq \frac {1}{16} \frac{\theta_2 \|x\|^p}{1-2^{p-3}}~,$$  for all
$x,y\in A$.
\end{corollary}
\begin{proof}
In Theorem 2.1, let $\varphi_1(x,y)=\theta_1 $ and
$\varphi_2(x,y)=\theta_2( \|x\|^p+\|y\|^p)$ for all $x,y \in A.$
\end{proof}

\begin{corollary}
Let $\theta_1$ and $\theta_2$ be nonnegative real numbers. Suppose
that
$$\|f(xy)-f(x)f(y)\|\leq \theta_1 ~,$$
$$\|f(2x+y)+f(2x-y)-2f(x+y)-2f(x-y)-12f(x)\|\leq \theta_2~,
$$ for all $x,y\in A$. Then there exists a unique cubic
homomorphism $~T:A \longrightarrow A$ such that
$$\|T(x)-f(x)\|\leq \frac{\theta_2}{14}~,$$ for all
$x\in A$.
\end{corollary}
\begin{proof}
It follows from Corollary 2.2.
\end{proof}

\begin{corollary}
Let $p\in (-\infty,3)$ and let $\theta$ be a positive real number.
Suppose that
$$\lim_{n \to \infty}\frac{\varphi(2^nx,2^ny)}{2^{6n}}=0~,$$
for all $x,y \in A.$ Moreover, Suppose that
$$\|f(xy)-f(x)f(y)\|\leq \varphi(x,y),$$
and that
$$\|f(2x+y)+f(2x-y)-2f(x+y)-2f(x-y)-12f(x)\|\leq \theta \|y\|^p,\eqno(2.7)$$
for all $x,y \in A.$ Then $f$ is a cubic homomorphism.
\end{corollary}
\begin{proof} Letting $x=y=0$ in $(2.7)$, we get that $f(0)=0.$ So
by $y=0$, in $(2.7),$ we get $f(2x)=2^3 f(x)$ for all $x \in A.$ By
using induction we have $$f(2^nx)=2^{3n} f(x),\eqno(2.8)$$ for all
$x \in A$ and $n \in \Bbb N.$ On the other hand by Theorem 2.1, the
mapping $T:A \to A$ defined by $$T(x)=\lim_{n \to \infty}
\frac{f(2^nx)}{2^{3n}},$$ is a cubic homomorphism. Therefore it
follows from $(2.8)$ that $f=T.$ Hence it is a cubic homomorphism.
\end{proof}

\begin{corollary}
Let $p,q,\theta\geq 0$ and  $p+q<3$. Let
$$\lim_{n \to \infty}\frac{\varphi(2^nx,2^ny)}{2^{6n}}=0~,$$
for all $x,y \in A.$ Moreover, Suppose that
$$\|f(xy)-f(x)f(y)\|\leq \varphi(x,y),$$
and that
$$\|f(2x+y)+f(2x-y)-2f(x+y)-2f(x-y)-12f(x)\|\leq \theta \|x\|^q \|y\|^p,$$
for all $x,y \in A.$ Then $f$ is a cubic homomorphism.
\end{corollary}
\begin{proof}
If $q=0$, then by Corollary 2.4 we get the result. If $q\neq 0,$ it
follows from Theorem 2.1, by putting $\varphi_1(x,y)=\varphi(x,y) $
and $\varphi_2(x,y)=\theta( \|x\|^p\|y\|^p)$ for all $x,y \in A.$

\end{proof}

\begin{corollary}
Let $p\in (-\infty,3)$ and $\theta$ be a positive real number. Let
$$\lim_{n \to \infty}\frac{\theta 2^{np} \|y\|^p}{2^{6n}}=0,$$
for all $x,y \in A.$ Moreover, suppose that
$$\|f(xy)-f(x)f(y)\|\leq \theta \|y\|^p~,$$
and
$$\|f(2x+y)+f(2x-y)-2f(x+y)-2f(x-y)-12f(x)\| \leq \theta \|y\|^p~,$$
for all $x,y \in A.$ Then $f$ is a cubic homomorphism.
\end{corollary}
\begin{proof}
Let $\varphi(x,y)=\theta \|y\|^p.$ Then by Corollary 2.4, we get the
result.
\end{proof}

\begin{theorem}
Let
$$\|f(xy)-f(x)f(y)\| \leq \varphi_1(x,y)~,$$
and
$$\|f(2x+y)+f(2x-y)-2f(x+y)-2f(x-y)-12f(x)\|\leq
\varphi_2(x,y)~, \eqno(2.9)$$ for all $x,y\in A$. Assume that the
series
$$\Psi(x,y)=\sum_{i=1}^{\infty}2^{3i}\varphi_2(\frac{x}{2^i},\frac{y}{2^i})$$ converges
and that
$$\lim_{n\longrightarrow
\infty}2^{6n}\varphi_1(\frac{x}{2^n},\frac{y}{2^n})=0,$$ for all
$x,y\in A$. Then there exists a unique cubic homomorphism $~T:A
\longrightarrow A$ such that
$$\|T(x)-f(x)\|\leq \frac{1}{16}\Psi(x,0)~, \eqno(2.10)$$ for all
$x\in A$.
\end{theorem}

\begin{proof}
Setting $y=0$ in $(2.9)$ yields
$$\|2 f(2x)-2.2^3 f(x)\|\leq \varphi_2(x,0) ~.\eqno (2.11)$$
Replacing $x$ by $\frac{x}{2}$ in $(2.11)$ to get
$$\|f(x)- 2^3
f(\frac{x}{2})\| \leq \frac{1}{2}\varphi_2(\frac{x}{2},0)~, \eqno
(2.12) $$ for all $x\in A$. By (2.12) we use iterative methods and
induction on $n$ to prove our next relation
$$\|f(x)-2^{3n}f(\frac{x}{2^n})\|\leq \frac{1}{2.2^3} \sum_{i=1}^{n} 2^{3i}
\varphi_2(\frac{x}{2^i},0)~.\eqno(2.13)$$ In order to show that the
functions $T_n(x)=2^{3n}f(\frac{x}{2^n})$ is a convergent sequence,
replace $x$ by $\frac{x}{2^m}$ in (2.13), and then multiplying  by
$2^{3m}$, where $m$ is an arbitrary positive integer. We find that
\begin{align*}
\|2^{3m}f(\frac{x}{2^m})-2^{3(n+m)} f(\frac{x}{2^{n+m}})\|&
 \leq \frac{1}{2.2^3}\sum_{i=1}^{n} 2^{3(i+m)}\varphi_2(\frac{x}{2^{i+m}},0)\\
&=\frac{1}{2.2^3}\sum_{i=1+m}^{n+m}
2^{3i}\varphi_2(\frac{x}{2^i},0)
\end{align*}
for all positive integers. Hence by the Cauchy criterion the limit
$~T(x)=\lim_{n\longrightarrow \infty}T_n(x)$ exists for each $x\in
A$. By taking the limit as $n\longrightarrow\infty$ in $(2.13)$, we
see that $\|T(x)-f(x)\|\leq \frac{1}{2.2^3}\sum_{i=1}^{\infty}
2^{3i} \varphi_2(\frac{x}{2^i},0)=\frac{1}{16}\Psi(x,0)$ and
$(2.10)$ holds for all $x\in A$. The rest of proof is similar to the
proof of Theorem 2.1.
\end{proof}

\begin{corollary}
Let $p>3$ and $\theta$ be a positive real number. Let
$$\lim_{n \to \infty} 2^{6n} \varphi(\frac{x}{2^n},\frac{y}{2^n})=0~,$$
for all $x,y \in A.$ Moreover, Suppose that
$$\|f(xy)-f(x)f(y)\|\leq \varphi(x,y)~,$$
and
$$\|f(2x+y)+f(2x-y)-2f(x+y)-2f(x-y)-12f(x)\|\leq \theta \|y\|^p~,\eqno(2.14)$$
for all $x,y \in A.$ Then $f$ is a cubic homomorphism.
\end{corollary}
\begin{proof} Letting $x=y=0$ in $(2.14)$, we get that $f(0)=0.$ So
by $y=0$, in $(2.14),$ we get $f(2x)=2^3 f(x)$ for all $x \in A.$ By
using induction we have
$$f(x)=2^{3n}f(\frac{x}{2^n}),\eqno(2.15)$$ for all $x \in A$ and $n
\in \Bbb N.$ On the other hand by Theorem 2.8, the mapping $T:A \to
A$ defined by $$T(x)=\lim_{n \to \infty}2^{3n}f(\frac{x}{2^n}),$$ is
a cubic homomorphism. Therefore it follows from $(2.15)$ that $f=T.$
Hence $f$ a cubic homomorphism.
\end{proof}

\paragraph{\large \bf Example}
Let
\[{\mathcal A} := \left[ \begin{array}{cccc}
{0} & {\Bbb R} & {\Bbb R} & {\Bbb R}\\
{0} & {0} & {\Bbb R} & {\Bbb R}\\
{0} & {0} & {0} & {\Bbb R}\\
{0} & {0} & {0} & {0}\\
 \end{array} \right], \]
then $\mathcal A$ is a Banach algebra equipped with the usual
matrix-like operations and the following norm:
\[ \|\left[ \begin{array}{cccc}
{0} & {a_1} & {a_2} & {a_3} \\
{0} & {0} & {a_4} & {a_5} \\
{0} & {0} & {0} & {a_6} \\
{0} & {0} & {0} & {0} \\
 \end{array} \right]\|=\sum_{i=1}^6| a_i | \hspace {00.1 cm} (a_i \in \Bbb R). \]
Let \[{a} := \left[
\begin{array}{cccc}
{0} & {0} & {1} & {2}\\
{0} & {0} & {0} & {1}\\
{0} & {0} & {0} & {0}\\
{0} & {0} & {0} & {0}\\
 \end{array} \right]\] and we define $f:\mathcal A \to \mathcal A$ by
$f(x)=x^3+a,$ and $$\varphi_1(x,y):=\|f(xy)-f(x)f(y)\|=\|a\|=4~,$$
$$\varphi_2(x,y):=\|f(2x+y)+f(2x-y)-2f(x+y)-2f(x-y)-12f(x)\|=14\|a\|=56~,$$
for all $x,y\in \mathcal A.$ Then we have
$$\sum_{k=0}^{\infty}\frac{\varphi_2(2^kx,2^ky)}{2^{3k}}=\sum_{k=0}^{\infty}\frac{56}{2^{3k}}=64~,$$
and
$$\lim_{n \to \infty}\frac{\varphi_1(2^nx,2^ny)}{2^{6n}}=0~.$$
Thus the limit $~T(x)=\lim_{n \to
\infty}\frac{f(2^nx)}{2^{3n}}=x^3$ exists. Also,
$$T(xy)=(xy)^3=x^3y^3=T(x)T(y)~.$$  Furthermore,
\begin{align*}
&T(2x+y)+T(2x-y)=(2x+y)^3+(2x-y)^3=16 x^3+12xy^2\\
&=2T(x+y)+2T(x-y)+12T(x)~.
\end{align*}
Hence  $T$ is cubic homomorphism.\\
Also from this example it is clear that the superstability of the
system of functional equations$$\hspace {1.5 cm}\begin{cases}
f(xy)=f(x)f(y), & \\
&\qquad\qquad.\\
f(2x+y)+f(2x-y)=2f(x+y)+2f(x-y)+12f(x), &
\end{cases}$$  with the control functions in
Corollaries 2.4, 2.5 and  2.6  does not hold.

\paragraph{\bf Acknowledgement.}
The authors would like to thank the referees for their valuable
suggestions. Also, the second  author would like to thank the office
of gifted students at Semnan University for its financial support.



\begin{thebibliography}{99}
\bibitem{Ao} T. Aoki, On the stability of the linear transformation in Banach spaces,
{\it J. Math. Soc. Japan.} 2(1950), 64-66.


\bibitem{BAD1} R. Badora, On approximate ring homomorphisms,{\it J. Math. Anal.
Appl.} \textbf{276}, 589-597 (2002).


\bibitem{B-L-Z} J. Baker, J. Lawrence and F. Zorzitto, The stability of the
equation $f(x+y)=f(x)f(y)$,{\it Proc. Amer. Math. Soc.}
\textbf{74} (1979), no. 2, 242-246.


\bibitem{Bo} D. G. Borugin, Class of transformations and bordering transformations,
{\it Bull. Amer. Math. Soc.} 27 (1951) 223-237.


\bibitem{F-R-s} V. A. Faizev, Th. M. Rassias and P. K. Sahoo, The space of
$(\psi,\gamma)$-additive mappings on semigroups,{\it Transactions of
the Amer. Math. Soc.} 354(11)(2002),4455-4472.


\bibitem{f.1} G. L. Forti, An existence and stability Theorem for a class
of functional equations,{\it Stochastica,} 4 (1980) 23-30.


\bibitem{f.2} G. L. Forti, Comments on the core of the direct method
for proving Hyers-Ulam stability of functional equations, {\it J.
Math. Anal. Appl,} 295 (2004), 127-133.



\bibitem{Hy} D. H. Hyers, On the stability of the linear functional equation,
{\it Proc. Natl. Acad. Sci. U.S.A.} 27 (1941) 222-224.


\bibitem{Hy-I-r} D. H. Hyers, G. Isac and Th. M. Rassias, Stability of Functional
Equations in Several Variables,{\it Birkhauser, Boston, Basel,
Berlin,} 1998.


\bibitem{Hy. Ra} D. H. Hyers and Th.M. Rassias, Approximate homomorphisms,
{\it Aequationes Math,} 44 (1992), 125-153.


\bibitem{J-K}  K. W. Jun  and H. M. Kim,  The generalized
 Hyers-Ulam-Russias stability of a cubic functional equation, {\it J. Math. Anal. Appl.} 274, (2002), no. 2, 267--278.

\bibitem{I-R} G. Isac and Th. M. Rassias, On the Hyers--Ulam stability of ?-additive
mappings, {\it J. Approx. Theory} 72(1993), 131-137.


\bibitem{Mal} L. Maligranda, A result of Tosio Aoki about a generalization of
Hyers-Ulam stability of additive functions- a question of
priority,{\it Aequat. Math.} 75 (2008) 289-296.



\bibitem{PAR1} C. Park, Hyers--Ulam--Rassias stability of homomorphisms in
quasi-Banach algebras,{\it Bull. Sci. Math.} \textbf{132} (2008),
no. 2, 87-96.



\bibitem{RAS1} Th. M. Rassias, The problem of S. M. Ulam for
approximately multiplicative mappings,{\it J. Math. Anal. Appl.}
{246} (2000), no. 2, 352-378.


\bibitem{Ra.t1} Th. M. Rassias and J. Tabor, Stability of mappings of Hyers-Ulam
Type,{\it Hadronic Press Inc. Florida,} 1994.


\bibitem{Ra.2} Th. M. Rassias, On a modified Hyers-Ulam sequence,
{\it J. Math. Anal. Appl.} 158(1991),106-113.


\bibitem{Ra.3} Th. M. Rassias, On the stability of functional equations originated
by a problem of Ulam,{\it Mathematica,} 44(67)(1)(2002),39-75.


\bibitem{Ra.4} Th. M. Rassias, On the stability of functional equations and a
problem of Ulam,{\it Acta Applicandae Math.} 62(1)(2000),23-130.


\bibitem {Ra.5} Th. M. Rassias, On the stability of the linear mapping in Banach
spaces,{\it Proc. Amer. Math. Soc.} 72 (1978) 297-300.



\bibitem {U} S. M. Ulam, A collection of Mathematical problems,
{\it Interscience Publ, New York,} 1960.
\end{thebibliography}
\end{document}